\documentclass{article}
\usepackage{arxiv}
\usepackage[utf8]{inputenc}
\usepackage[T1]{fontenc}
\usepackage{url}
\usepackage{booktabs}
\usepackage{amsfonts}
\usepackage{nicefrac}
\usepackage{microtype}
\usepackage{graphicx}
\graphicspath{{./images/}}
\usepackage{xcolor}
\usepackage{cite}
\usepackage{mathtools}
\usepackage{amssymb}
\usepackage{amsthm}
\usepackage{hyperref}

\newtheorem{theorem}{Theorem}
\newtheorem{proposition}{Proposition}
\theoremstyle{definition}
\newtheorem{definition}{Definition}
\newtheorem{example}{Example}
\newenvironment{credits}{\endgraf\begingroup\small}{\endgraf\endgroup}
\newcommand{\ackname}{Acknowledgments.}

\newcommand{\inst}[1]{\textsuperscript{#1}}
\newcommand{\orcidID}[1]{\,{\normalfont\scriptsize\href{https://orcid.org/#1}{(#1)}}}
\newcommand{\email}[1]{{\normalfont\ttfamily #1}}
\newcommand{\mciinstitutes}{}
\newcommand{\institute}[1]{\gdef\mciinstitutes{#1}}
\newcounter{mciinstitute}
\let\arxivauthor\author
\renewcommand{\author}[1]{\arxivauthor{%
  \begin{minipage}{0.96\textwidth}\centering
    {\def\and{\\[2pt]}#1\endgraf}
    \medskip
    {\normalfont\small
      \setcounter{mciinstitute}{1}%
      \def\and{\endgraf\medskip\stepcounter{mciinstitute}\textsuperscript{\arabic{mciinstitute}}}%
      \textsuperscript{1}\mciinstitutes\endgraf}
  \end{minipage}}}
\newcommand{\titlerunning}[1]{\renewcommand{\shorttitle}{#1}}
\newcommand{\authorrunning}[1]{\renewcommand{\headeright}{#1}}
\date{}

    \newcommand{\Cl}{\mathcal{C}l}
    \newcommand{\Op}{\mathcal{O}p}
    \newcommand{\ep}{\varepsilon}
    \newcommand{\Bt}[2]{(B_{#1})_{#2}}
    \newcommand{\R}{\mathbb{R}}

    \DeclareMathOperator{\MCI}{MCI}
    \DeclareMathOperator{\conv}{conv}
    \DeclareMathOperator{\inter}{int}

\begin{document}
    \title{Multiscale Closing Index (MCI) for measuring structural similarity between 2D point sets} 
    %
    %
    \author{Nicol\`o Zava\inst{1}\orcidID{0000-0001-8686-1888} \and Angel Franco-Alvarez\inst{2}
    \and Samy Blusseau\inst{3}\orcidID{0000-0003-0294-8172}   \and
    Santiago Velasco-Forero\inst{3}\orcidID{0000-0002-2438-1747}\and 
    Maria-Jose Jimenez\inst{2}\orcidID{0000-0001-7316-117X}
    }
    \authorrunning{Zava et al.}
    \titlerunning{Multiscale closing index}
    \institute{
    Institute of Science and Technology Austria (ISTA)\\
    \email{nicolo.zava@gmail.com}\and
    Universidad de Sevilla\\ 
    \email{afalvarez95@gmail.com, \,majiro@us.es}\and
    Mines Paris, PSL University, Centre for Statistics and Images
    (STIM), Fontainebleau, France\\ 
    \email{\{samy.blusseau, santiago.velasco\}@minesparis.psl.eu}    
    }
    \maketitle   
    %
    \begin{abstract}
    We introduce Multi-scale Closing Index (MCI), a new dissimilarity measure for evaluating the alignment of two 2D point sets. MCI is bounded between 0 and 1 and is based on multiscale mathematical morphology operations, with smaller values indicating greater structural agreement. The measure was motivated by the need to assess the quality of registration between point sets representing cell populations in biological images. 
    To evaluate its behaviour under controlled conditions, we generated synthetic point clouds from Perlin-noise-based spatial structures, following established approaches for simulating biologically plausible spatial organization. We assessed MCI between those point clouds and different transformed versions, obtained by local perturbation, point removal, translation, rotation and addition of impulsive noise. Across these experiments, the dissimilarity MCI increased gradually and consistently with perturbation severity, while showing a clear discontinuity between clouds originating from the same underlying spatial structure and point clouds generated from different structures. We further compared MCI with commonly used point-set distances, including Hausdorff-based measures and Sliced Wasserstein distance. These experiments show the ability of MCI to capture structural agreement and quantify deviations between point sets, highlighting its potential as a complementary metric for evaluating two-dimensional point-set registration.
    
    \keywords{2D point sets similarity measure \and alignment quality score\and morphological closing \and Jaccard index.}
    \end{abstract}

\section{Introduction}

Point-set registration aims to align two or more point sets by estimating correspondences and a spatial transformation that brings them into agreement \cite{Zhu2019Review,Maiseli2017Trends,Pei2023Review}. It is widely used in computer vision, robotics, medical imaging, and spatial data analysis, but remains challenging in the presence of noise, outliers, missing points, partial overlap, and heterogeneous sampling.

Evaluating the quality of an alignment is itself nontrivial because different measures capture different notions of disagreement. Nearest-neighbour distances quantify local geometric error, whereas distributional distances characterize the displacement of the point sets more globally. More elaborate kernel, probabilistic, graph-based, or correspondence-based criteria also exist, but many are designed as components of the registration procedure itself rather than as simple post-registration similarity measures. Here, we focus on three complementary references for already aligned 2D point sets.

We consider the 95th-percentile Hausdorff distance (HD95), which captures large local discrepancies while reducing the strong outlier sensitivity of the classical Hausdorff distance \cite{Eiter1997Distance}; the Average Hausdorff Distance (AHD), which summarizes typical nearest-neighbour disagreement across the sets \cite{Aydin2021Hausdorff}; and the Sliced Wasserstein distance, which compares the two point sets as empirical spatial distributions and therefore provides a more global measure of displacement.

However, these distances do not explicitly quantify whether the spatial organization represented by the point clouds is preserved across different length scales. This is particularly relevant in spatial biology, where alignment should preserve meaningful organization beyond pointwise proximity \cite{Bull2024Colorectal}. This motivates assessing alignment through a representation that progressively incorporates spatial context, from fine local structure to broader organization.

We therefore introduce a multiscale morphological measure for comparing aligned 2D point sets. The proposed measure applies morphological operations over increasing radii and integrates the resulting agreement into a single dissimilarity score, thereby assessing how similarly the structures induced by the two point sets evolve across spatial scales. We compare it with HD95, AHD, and Sliced Wasserstein, selected as compact and complementary benchmarks for local extreme, local average, and global distributional disagreement.

In this paper, Section~\ref{sec:background} reviews similarity measures for evaluating point-set alignment, as well as basics of mathematical morphology operators. Section~\ref{sec:method} introduces the proposed multiscale morphological structural measure. Section~\ref{sec:experiments} presents the experimental evaluation and comparative analysis. Section~\ref{sec:conclusions} discusses the strengths and limitations of the approach, and concludes the paper.

\section{Background}\label{sec:background}

In this section, first, we introduce several methods used in the literature for measuring 2D point set similarity. Secondly, we describe
the basic tools of mathematical morphology that will be used in the rest of the paper.

\paragraph{Similarity measures for 2D point sets.}
\label{sec:alignment_similarity}

Let $X=\{x_i\}_{i=1}^{n}$, $
Y=\{y_j\}_{j=1}^{m}$
be two finite point sets in $\R^2$. We consider three established dissimilarity  measures, with smaller values indicating better alignment: percentile Hausdorff, average Hausdorff, the Sliced Wasserstein distances. We write $d(x,Y)=\min_{y\in Y}\|x-y\|_2$ for the distance from a point to a set.

\textit{Hausdorff and percentile Hausdorff distances.}
The symmetric Hausdorff distance
\begin{equation*}
d_H(X,Y)
=
\max
\big\{
\max_{x\in X} d(x,Y),
\max_{y\in Y} d(y,X)
\big\}
\end{equation*}
is a classical metric for finite point sets
\cite{Huttenlocher1993Hausdorff,Eiter1997Distance}, but, being determined by the single largest nearest-neighbour discrepancy, it is highly sensitive to outliers, missing points, and sampling artefacts \cite{Eiter1997Distance,Taha2015Metrics}. We therefore use the 95th-percentile Hausdorff distance (\textit{HD95}). Let $\mathcal{D}_{X\rightarrow Y}
=
\left\{
d(x,Y)\mid x\in X
\right\}
$,
and define $\mathcal{D}_{Y\rightarrow X}$ analogously. Then
\begin{equation*}
d_{\mathrm{HD95}}(X,Y)
=
\max
\left\{
P_{95}\!\left(\mathcal{D}_{X\rightarrow Y}\right),
P_{95}\!\left(\mathcal{D}_{Y\rightarrow X}\right)
\right\},
\end{equation*}
where $P_{95}$ denotes the empirical 95th percentile. HD95 retains sensitivity to large misalignments, while reducing the influence of a few extreme points, and it is a standard robust alternative in segmentation evaluation \cite{Taha2015Metrics}. 

\textit{Average Hausdorff distance (AHD).}
It captures the overall nearest-neighbour discrepancy, rather than its upper tail:
\begin{equation*}
d_{\mathrm{AHD}}(X,Y)
=
\frac{1}{2}
\bigg(
\frac{1}{n}
\sum_{i=1}^{n} d(x_i,Y)
+
\frac{1}{m}
\sum_{j=1}^{m} d(y_j,X)
\bigg).
\end{equation*}
It is considerably less sensitive to individual extreme points than $d_H$ \cite{Dubuisson1994Hausdorff,Aydin2021Hausdorff}.

\textit{Sliced Wasserstein distance.}
The Sliced Wasserstein distance provides an efficient approximation to optimal-transport comparison by projecting the two point sets onto a set of one-dimensional directions and computing the Wasserstein distance between the corresponding projected distributions \cite{Bonneel2015Sliced,Peyre2019ComputationalOT,Villani2003Topics,Santambrogio2015OT}. For a unit direction $\theta$, let
\[
X_{\theta}=\{\langle\theta,x_i\rangle\}_{i=1}^{n},
\qquad
Y_{\theta}=\{\langle\theta,y_j\rangle\}_{j=1}^{m}
\]
denote the projected point sets. Using $L$ directions $\theta_1,\ldots,\theta_L$, we approximate the 2-Sliced Wasserstein distance by
\begin{equation*}
d_{\mathrm{SW}}(X,Y)
=
\bigg(
\frac{1}{L}
\sum_{\ell=1}^{L}
W_2^2(X_{\theta_\ell},Y_{\theta_\ell})
\bigg)^{1/2},
\end{equation*}
where $W_2$ denotes the one-dimensional 2-Wasserstein distance between the empirical distributions of the projected points. Thus, unlike nearest-neighbour measures, Sliced Wasserstein captures differences in the overall spatial distribution of the two point sets.

None of the above measures explicitly evaluates agreement over a sequence of spatial scales. This motivates the proposed multiscale morphological criterion, which compares the point sets across progressively increasing morphological scales and aggregates these contributions into a single alignment score.


\paragraph{Basic concepts of MM.}
Mathematical morphology is a branch of mathematics which combines elements of fields such as set theory, lattice algebra and integral geometry to analyze geometrical structures \cite{Soille2010Morphological} especially in image analysis \cite{najman2013mathematical}. 
Mathematical morphology studies operations that modify the shape of sets in a mathematical space, with dilation and erosion being two of its fundamental operations (see \cite{Matheron1975,Serra1982} for foundational references).

Let \(X\subseteq \R^{2}\) be a set and let $B_r$ be the closed ball of radius ($r\geq 0$) in $\R^2$, centered at the origin.

For any \(x\in  \mathbb{R}^{2}\), the translation of \(B_r\) by \(x\) is denoted by $\Bt{r}{x}\coloneqq\{b+x\mid b\in B_r\}$.
The \emph{dilation} and \emph{erosion} of \(X\) by \(B_r\), denoted respectively $\delta_r(X)$ and $\ep_r(X)$, are defined by \cite{Soille2010Morphological}:
\begin{center}
\begin{minipage}[b]{.45\textwidth}
\vspace{-\baselineskip}
    \begin{equation*}
        \delta_r(X) = \{x \mid \Bt{r}{x} \cap X \ne \emptyset\}
\end{equation*}
\end{minipage}
 \hfill and \hfill
\begin{minipage}[b]{.4\textwidth}
\vspace{-\baselineskip}
\begin{equation*}
    \ep_r(X) = \{x \mid \Bt{r}{x} \subseteq X\}.
\end{equation*}
\end{minipage}
\end{center}

Thus, the dilation of $X$ by $B_r$ is the set of all points whose distance to $X$ is at most $r$, whereas the erosion of $X$ by $B_r$ is the set of all points $p$ such that the closed ball $(B_r)_p$ of radius $r$ centered at $p$ is entirely contained in $X$. Furthermore, $\delta_r(X) = \ep_r(X^C)^C$ where $A^C \coloneqq \mathbb{R}^2\setminus A$ denotes the complement of $A.$

The composition of $\delta_r$ followed by $\ep_r$ defines another fundamental morphological operation, called the \emph{closing by the structuring element $B_r$}:
\begin{equation*}
    \Cl_{r}(X) = \ep_r(\delta_{r}(X)).
\end{equation*}

The closing operator has the following fundamental properties. For all $X,Y\subseteq\mathbb{R}^2$ and $r\geq 0$:
\begin{enumerate}
\item Extensivity:
$X\subseteq \Cl_r(X).$
Thus, closing can fill gaps and bridge narrow connections, but does not remove points from $X$.
\item Increasingness with respect to sets:
$X\subseteq Y \,\Longrightarrow\,\Cl_r(X)\subseteq\Cl_r(Y).$
Hence, the inclusion relation between sets is preserved by closing.
\item Increasingness with respect to the radius: $ \forall s\geq r, \ \Cl_r(X) \subseteq \Cl_s(X)$. 
\item Idempotence:
$\Cl_r\bigl(\Cl_r(X)\bigr)=\Cl_r(X).$
Thus, applying the same closing operation more than once does not produce any further change.
\item Closing and convexity \cite{Matheron1975}: if $X$ is closed then, for any $r\geq 0$, $\Cl_r(X) \subseteq \conv(X)$ where $\conv(X)$ is the convex hull of $X$. Hence if $X$ is also convex, then $\Cl_r(X) = X \ \ \forall r\geq 0$ .
\end{enumerate}
Furthermore, $\Cl_r$ is the dual of the opening $\Op_r$, defined by $\Op_r(X)\coloneqq \delta_r(\ep_r(X))$, in the sense that $\Cl_r(X) = \Op(X^C)^C$. Since $\Op_r(X) = \bigcup
_{x\in\ep_r(X)}\Bt{r}{x}$, we get the following geometric interpretation for the closing:
\begin{equation*}
    \Cl_r(X) = \bigcap_{x\notin\delta_r(X)}\Bt{r}{x}^C = \bigcap\{\Bt{r}{x}^C \mid \Bt{r}{x}\cap X=\emptyset\}.
\end{equation*}

In words, $\Cl_r(X)$ is the intersection of all the complement of closed balls which do not intersect $X$. Note that, for $r>0$, this is also what is called the $\alpha$-hull of $X$ \cite{refId0}, for $\alpha = -1/r$.

\section{Multiscale Closing Index}\label{sec:method}

In this section, we introduce our dissimilarity measure, Multiscale Closing Index (MCI), together with some theoretical results. 

Given two bounded measurable subsets $X$ and $Y$ of $\mathbb R^2$ with at least one of them with non-zero volume, their {\em Jaccard index} is the value
\[
J(X,Y)=\frac{\lambda(X\cap Y)}{\lambda(X\cup Y)},
\]
where $\lambda$ is the Lebesgue measure on the plane (whenever $\lambda(X\cup Y)>0$). Clearly, $0\leq J(X,Y)\leq 1.$ For every pair of bounded measurable subsets $X$ and $Y$ of the plane, we define their {\em Jaccard distance} as
\[
J_\delta(X,Y)=
\begin{cases}
\begin{aligned} 
& 1-J(X,Y) &\text{if $\lambda(X\cup Y)\neq 0$},\\
& 0  
&\text{otherwise.}
\end{aligned}
\end{cases}
\]

\begin{theorem}[\cite{Marczewski1958}]\label{thm:Jaccard_distance_is_distance}
On the family of measure bounded subsets of the plane up to subsets with $0$-measure, $J_\delta$ is a metric. More explicitly, $J_\delta$ is symmetric, satisfies the triangle inequality, and, for any pair of subsets $X$ and $Y$ of the plane, $J_\delta(X,Y)=0$ if and only if their \textit{symmetric difference} $X\triangle Y=(X\cup Y)\setminus(X\cap Y)$ has measure $0$.
\end{theorem}

In our setting, we are interested in applying the Jaccard metric to compare point clouds at different scales. First, given two point clouds $X$ and $Y$ in the plane, and a certain scale $r>0$, we define two subsets $X_r$ and $Y_r$ and compare them using the Jaccard distance. We then let the scale vary and summarise the family $\{J_\delta(X_r,Y_r)\}_r$ in one index. In this way, we preserve some crucial properties of the Jaccard distance, but keeping into account the behaviour at different scales. As for the first step, we will use the closing as defined above. The following result will be useful to analyze some properties of our index.

\begin{proposition}\label{prop:closing_to_conv}
    Let $X\subset\R^2$ be compact, and let $\conv(X)$ denote its convex hull. Then
    \begin{equation*}
    \lim\limits_{r\to+\infty}J_\delta\big(\Cl_r(X),\conv(X)\big) = 0.
    \end{equation*}
    If furthermore $X$ contains at least three non-colinear points or is convex, then the convergence also holds for the Hausdorff distance
        \begin{equation*}
    \lim\limits_{r\to+\infty} d_H \big(\Cl_r(X),\conv(X)\big) =  0.
    \end{equation*}
\end{proposition}
Results on the convergence of closings to the convex hull have been known for long \cite{Matheron1975,Serra1982} but they were stated in the hit-or-miss topology, for which, like for the Hausdorff one, the non-convex sets of aligned points are counter-examples (e.g. if $X$ is a pair of points, $\Cl_r(X) = X \neq \conv(X)$ for any $r\geq 0$). Hence, the motivation of this proposition is to highlight this case and to state that it is not an issue with the Jaccard distance.

For a compact set $X$ containing non-colinear points, the result comes from the fact that 1. its convex hull has non-empty interior and 2. that any point $x \in \inter(\conv(X))$ is included in a closing $\Cl_r(X)$ for $r$ large enough. This is obvious for a finite set from the geometric interpretation of closings, and it can be generalized to any compact set thanks to Steinitz's theorem \cite{Steinitz1913}. The case of a convex compact set $X$ of aligned points is obvious because then $\Cl_r(X) = X = \conv(X)$. The convergence holds with the Jaccard distance for a set of non-convex aligned points because any of its closings, as well as its convex hull, have zero Lebesgue measure.

We can now move on to the definition of our index. Let $X$ and $Y$ be two point clouds in the plane and a fixed radius $r\geq 0$. 
The value $J_\delta(\Cl_r(X),\Cl_r(Y))$ is still scale-dependent as it looks only at a specific scale. We assemble these values by taking their normalised integral. 

\begin{definition}
    Let $X$ and $Y$ be two point clouds in $\mathbb R^2$. The {\em Multiscale Closing Index of $X$ and $Y$ from radius $r$ up to radius $R$}, where $R\geq r\geq 0$, is defined as
    \[
    \MCI_r^R(X,Y)=\frac{1}{R-r}\int_{r}^RJ_\delta(\Cl_s(X),\Cl_s(Y))ds.
    \]
    If $r=0$, we simply write $\MCI^R(X,Y)$.
\end{definition}

\begin{figure}[htbp]
    \centering
    \includegraphics[width=1\textwidth]{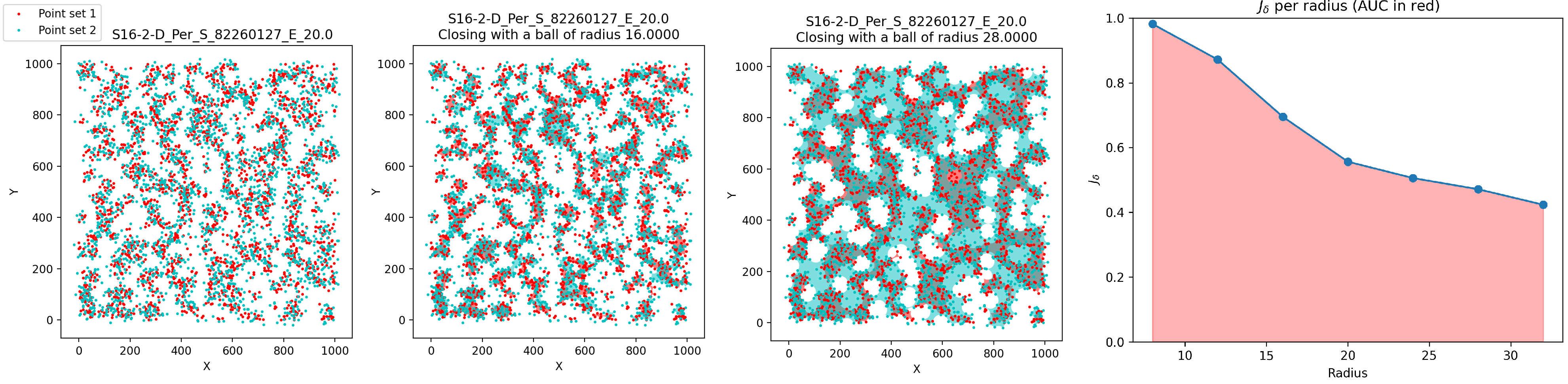}
    \caption{From left to right: two point clouds $X$ and $Y$; overlapping of subsequent closings with increasing radius, 16 and 28; graph of values $J_\delta(\Cl_r(X),\Cl_r(Y))$ for the list of radii $\{8,12,16,20,24,28,32\}$.}
    \label{fig:closings}
\end{figure}

\begin{example}
    In Fig. \ref{fig:closings}, two sets of points are used to illustrate how the \emph{shapes} change when the radius of the closing is increased. Intuitively, the MCI ($0.63$) is the proportion of the area of the rectangle of width $(32-8)$ and height 1, that is under the curve.
\end{example}

The Multiscale Closing Index satisfies the following properties.

\begin{proposition}\label{prop:MCI_metric}
    Let $R> r\geq 0$. 
    \begin{enumerate}
        \item $\MCI_r^R$ takes values between $0$ and $1$. 
        \item $\MCI_r^R$ is symmetric and satisfies the triangle inquality.
        \item Given two point clouds $X$ and $Y$, $\MCI_r^R(X,Y)=0$ if there is a finite set $I\subseteq[r,R]$ such that, for every $s\in [r,R]\setminus I$
        , $\lambda(\Cl_s(X)\triangle \Cl_s(Y))=0$.
    \end{enumerate}
\end{proposition}
\begin{proof}
    Item 1. and the symmetry part of item 2. immediately follows from the definitions. Let us show the triangle inequality. Theorem \ref{thm:Jaccard_distance_is_distance} states, in particular, that $J_\delta$ satisfies the triangle inequality, and the pointwise inequality is preserved once the integral and the normalisation are applied.

    Item 3. follows from the fact that the integrand function $J_\delta(\Cl_s(X),\Cl_s(Y))$ is $0$ if and only if $\lambda(\Cl_s(X)\triangle\Cl_s(Y))=0$.
\end{proof}

In the direction of item 3 of Proposition \ref{prop:MCI_metric}, the reader may wonder if the Multiscale Closing Index is able to distinguish points, i.e., being zero only if the two point clouds are the same. The $\MCI$, however, does not satisfy this property for two concurring reasons. First, already the Jaccard metric can distinguish points only up to zero-measure: if $X$ and $Y$ are two subsets whose symmetric difference has measure zero, then $J_\delta(X,Y)=0$  (Theorem \ref{thm:Jaccard_distance_is_distance}). Furthermore, the closings of two different point clouds may coincide as soon as they have positive measure, as the following example shows.
\begin{example}
    Let $X$ be the vertices of an equilateral triangle and $Y$ be $X$ together with the center of the triangle. Below a certain radius both closings have zero measure and, for larger radii the closings are equal. Hence for every $R\geq 0$, $\MCI^R(X,Y)=0$.
\end{example}

Let us present a few more properties of the Multiscale Closing Index in relation to the convex hull. We shall first notice that, for a compact set $X$ of the plane, containing three non-colinear points is equivalent to $\inter(\conv(X))\neq\emptyset$, and also equivalent to the existence of $r_0\geq 0$ such that $\lambda(\Cl_r(X))>0$ for any $r\geq r_0$. For such a set, which we will call \emph{generic}, we will note $r(X)\coloneqq \inf\lbrace r\geq 0\mid \lambda(\Cl_r(X))>0 \rbrace$.

\begin{proposition}
    Let $X,Y$ be two point clouds in $\mathbb R^2$.
    \begin{enumerate}
        \item For every $r\geq 0$, $\MCI_r^R(X,Y)\xrightarrow{R\to\infty}J_\delta(\conv(X),\conv(Y))$.
        \item If $X$ and $Y$ are generic, then: $\lambda(\conv(X)\cap\conv(Y))=0 \iff \MCI^R_r(X,Y)=1$ for every $R>r>r_0(X,Y)\coloneqq \max(r(X), r(Y))$.
    \end{enumerate}
\end{proposition}
\begin{proof}
    1. It follows from the fact that the closings of $X$ and $Y$ converge to the respective convex hulls (see Proposition \ref{prop:closing_to_conv}).
    
    2. If $\lambda(\conv(X)\cap\conv(Y))=0$, then, for every $r\geq 0$, $\lambda(\Cl_r(X)\cap\Cl_r(Y))=0$, and so, for $R\geq s\geq r> r_0(X,Y)$, $J_\delta(\Cl_s(X),\Cl_s(Y))=1$, which implies that $\MCI_r^R(X,Y)=1$. Conversely, if $\lambda(\conv(X)\cap\conv(Y))>0$, then $J_\delta(\conv(X),\conv(X))<1$ and, given Proposition \ref{prop:closing_to_conv}, 
    $\exists A\leq0$ s.t. $s\geq A \Rightarrow J_\delta(\Cl_s(X),\Cl_s(Y))<1$.
    Then for $R>r$ large enough, $\MCI_r^R(X,Y)<1$.
\end{proof}

\section{Experiments and results}\label{sec:experiments}

We designed a set of experiments to show the robustness and coherence of our index. 

\paragraph{Synthetic dataset.} We generated 48 synthetic two-dimensional point sets using Perlin noise
\cite{perlin1985image,perlin2002improving}, which produces spatially
correlated random fields and has been used to simulate biological and
tissue-like spatial structures
\cite{gharia2020signal,najafiaghdam20223d,LAWSON2024perlin,li2026synora}.
All point sets were generated on a $1000\times1000$ pixel domain using
spatial-frequency parameters $o\in\{4,8,16\}$, corresponding
approximately to characteristic spatial scales
$\lambda \approx 1000/o=\{250,125,62.5\}$ pixels.

\begin{figure}[htbp]
    \centering
    \includegraphics[width=0.7\textwidth]{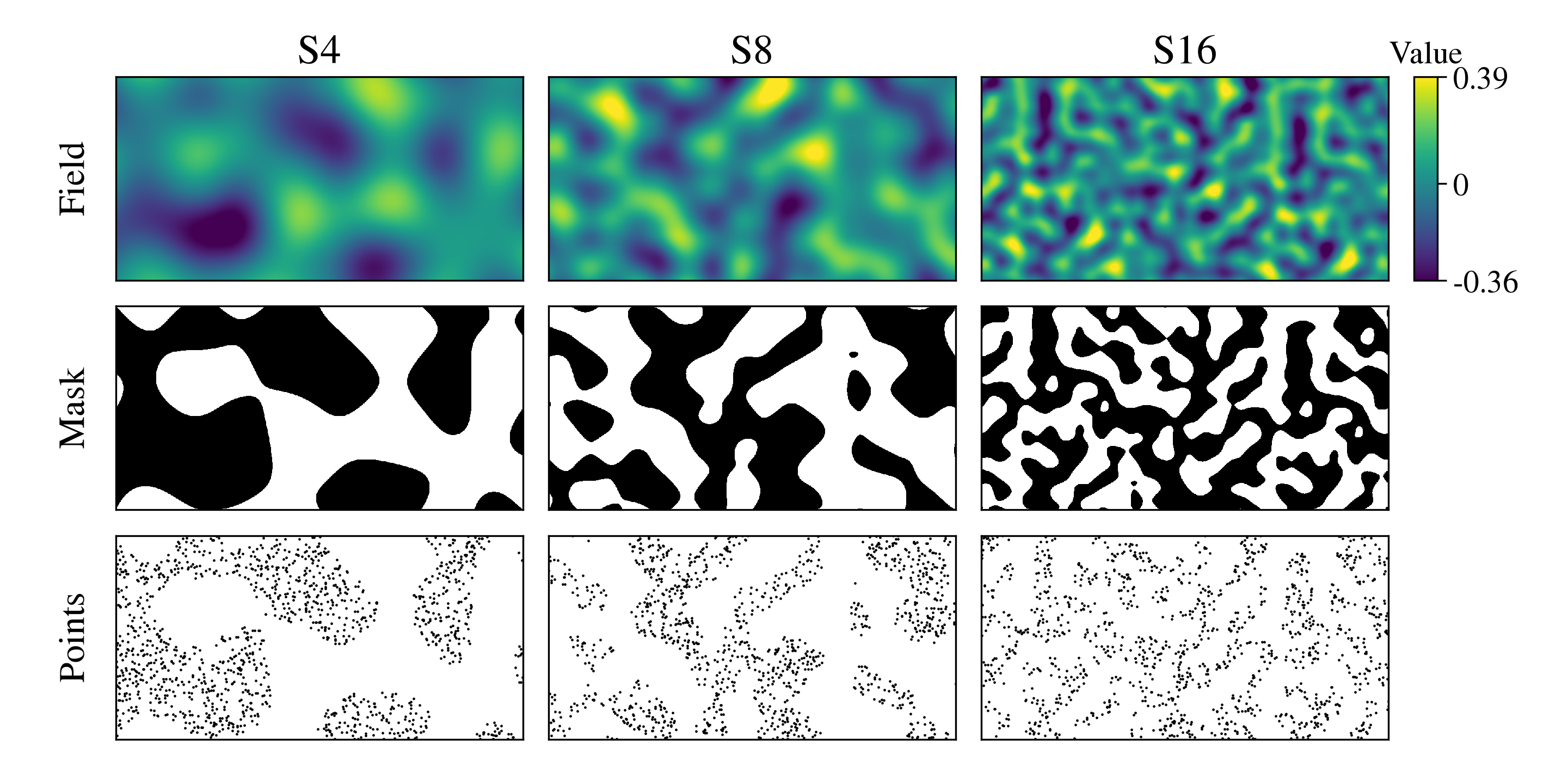}
    \caption{
Illustration of the synthetic point-cloud generation process using Perlin noise.
Columns show three representative values of the Perlin spatial-frequency parameter, $o=4,8$ and $16$, respectively. Increasing $o$ produces progressively finer spatial structure. First row: continuous Perlin fields; second row: their binary masks obtained by median thresholding; third row: resulting point clouds
sampled from the black pixels with fixed point density.
}    \label{fig:perlin}
\end{figure}

Twelve underlying spatial configurations were generated and each was sampled four times independently, yielding 48 point sets. Three types of
fields were considered:
\begin{itemize}
    \item \textit{Simple (S4, S8, S16):} two independent single-component fields for
    each $o\in\{4,8,16\}$ (see Fig. \ref{fig:perlin} for an example of each one);
    
    \item \textit{Multiweight Balanced (B4816):} two realizations of a  multiscale field combining components with $o=\{4,8,16\}$ and equal weights
    $(1/3,1/3,1/3)$;
    
    \item \textit{Multiweight Custom (C416-Coarse, C416-Fine):} two multiscale fields combining $o=4$ and $=16$, using coarse-dominated weights $(0.65,0.35)$ and fine-dominated weights $(0.35,0.65)$, with two independent realizations of each configuration.
\end{itemize}

Different seeds were used for each Perlin component. For multicomponent fields, components were combined as a weighted sum of the raw Perlin noise values, and the resulting field was not further normalized. Each field was binarized at its median value, yielding approximately $50\%$ candidate pixels. Four independent point sets were then obtained from each binary region by
sampling exactly $0.5\%$ of its candidate pixels uniformly without replacement, using different sampling seeds. This resulted in approximately $2500$ points per point set. Coordinates were stored in $(x,y)$ pixel coordinates.

This design provides controlled single and multiscale spatial structures while allowing direct comparison between independent point samples drawn from the same underlying region and samples generated from different regions.

\paragraph{Experiments design.} An  initial baseline experiment was conducted to verify that the proposed similarity measure does not yield spuriously high similarity scores for uncorrelated, random point clouds. The metric's performance was subsequently evaluated against the HD95, AHD and Sliced Wasserstein distances. 

Second, for each set $P_i$ in the collection of generated point sets $\mathcal{C}=\{P_i\}_{i\in I}$, $I=\{1,\ldots , 48\}$, we performed the following operations to produce a new point set $\phi(P_i)$. Then, we computed MCI for each pair $(P_i, \phi(P_i))$, where $\phi$ accounts for each of the following operations:
\begin{itemize}
    \item \textit{Perturbation 
    }: Points are generated by applying a perturbation on the coordinates. That is, fixing \(\varepsilon>0\), for each point \(x\), choose a random unit vector \(v\) in the plane and move \(x\) by distance \(\varepsilon\) in that direction:
    $x \mapsto x+\varepsilon v.$ We implemented this operation for $\varepsilon \in\{5,10,20,40\}$.
    \item \textit{Points removal:} Removal of a proportion $(10\%, 20\%, 40\% \mbox{ or } 60\%)$ of the initial amount of points is performed randomly. 
    \item \textit{Perturbation and removal:} Points are obtained by combining both perturbation and removal as explained before. 
    \item \textit{Translation
    :} Points are generated by applying a translation along the horizontal axis by $\varepsilon$ pixels with $\varepsilon \in \{5,10,20,40\}$.
    
    \item \textit{Rotation
    :} Points are generated by applying a rotation with respect to the central point of the domain with a fixed angle from the set $\vartheta\in\{2^\circ, 5^\circ, 7^\circ, 10^\circ\}$.
   \item \textit{Impulsive noise addition:} additional points (25, 50, 75, and 100) were added to each of the 48 base Perlin clouds, sampled from the white regions of the binary masks, while all original points were retained. 
\end{itemize}

Regarding the minimum and maximum radius, $r$ and $R$, for MCI computations, they were calculated based on the local sampling scales of the point clouds using their second-nearest-neighbor (2NN) distance distributions.  
The lower and upper radii were selected as the maxima, across clouds involved in each experiment, of the 5th and 95th percentiles, respectively. This provides a common interval accommodating differences in sampling density while reducing the influence of extreme local spacings. The use of 2NN distances is motivated by three-point neighborhoods, the smallest configurations that can produce nontrivial closings. These distances serve as practical scale estimates rather than exact thresholds for morphological transitions.

\begin{figure}[htbp]
    \centering
    \includegraphics[width=1\textwidth]{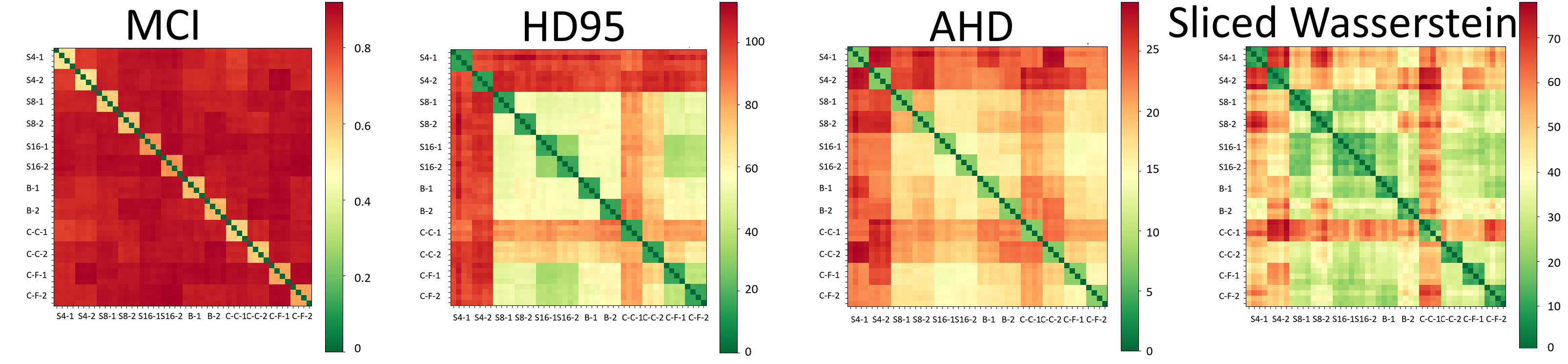}
    \caption{Pairwise dissimilarities among 48 Perlin clouds: MCI, HD95, AHD and sliced Wasserstein. Each diagonal block contains four independent samplings of one region.}
    \label{fig:heatmaps}
\end{figure}

\paragraph{Results.} Starting with our baseline experiment, we can visualize the cross-comparison of our 48 Perlin point clouds via heatmaps in Fig. \ref{fig:heatmaps}, which suggests that our measure provides a more robust similarity assessment than the alternative methods. In particular, MCI provides a consistent score interpretation across the tested spatial configurations: independent samplings of the same regions receive intermediate dissimilarities, whereas samplings of different regions receive uniformly high dissimilarities. The alternative measures exhibit greater variation among different-region comparisons, making a common interpretation of poor alignment less straightforward.

\begin{figure}[htbp]
    \centering
    \includegraphics[width=1\textwidth]{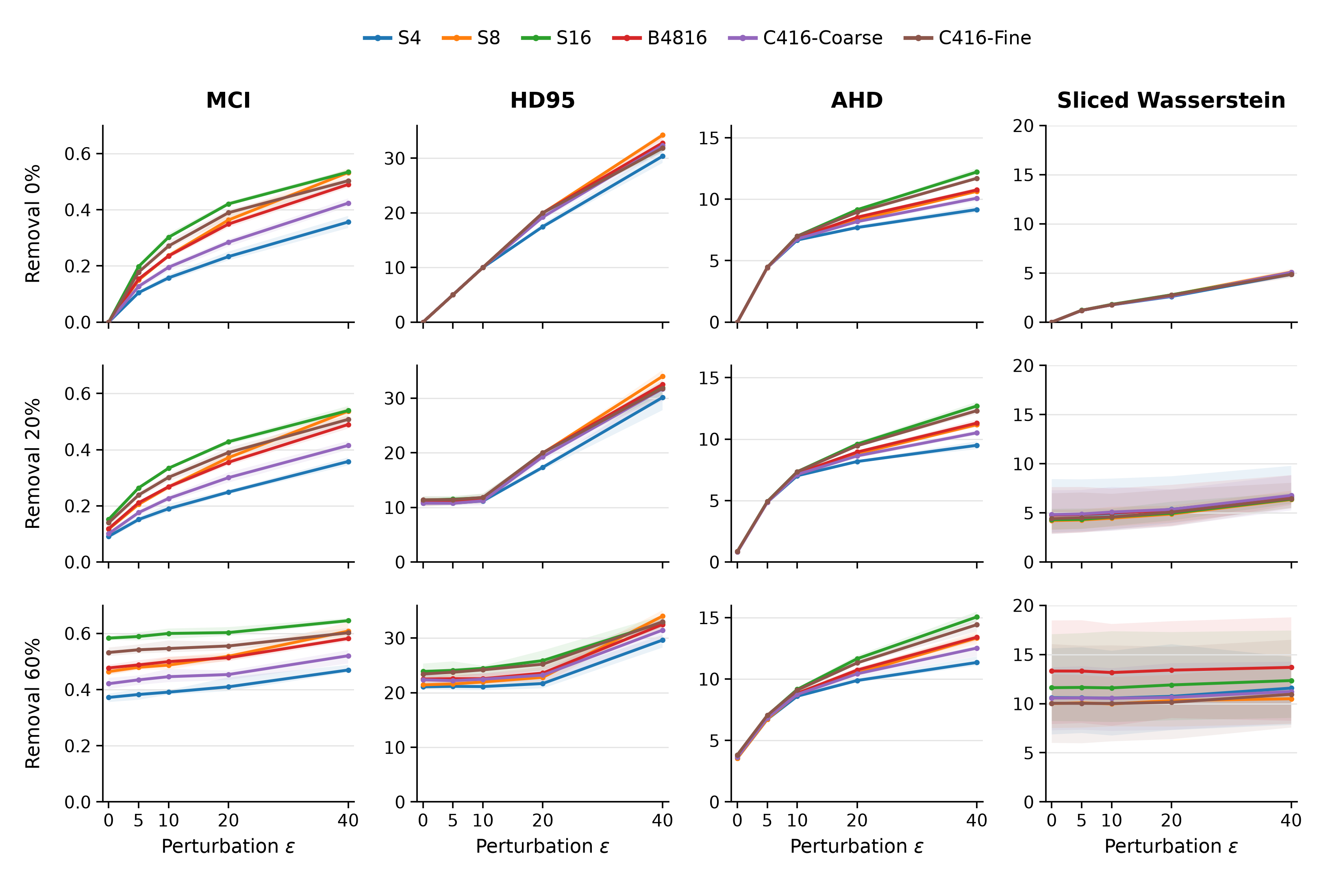}
    \caption{Perturbation response at 0\%, 20\% and 60\% removal. Colours identify cloud classes; lines show means and bands show observed ranges across eight clouds per class (two regions, four samplings each). Full removal levels are provided in the repository.}
    \label{fig:per+rem}
\end{figure}

Across the transformation experiments, MCI generally increases with severity while retaining clear sensitivity to structural scale: comparable alterations produce larger dissimilarities in finer structures. Under translation, sliced Wasserstein assigns identical scores to all configurations, whereas MCI distinguishes them throughout. At $10^\circ$, HD95 and sliced Wasserstein assign larger scores to S4 than S16. The combined experiment (Fig.~\ref{fig:per+rem}) reinforces MCI's coarse--fine differentiation, particularly at mild perturbations, where HD95 and AHD provide little configuration separation. Sliced Wasserstein shows greater variation between realizations under removal and occasionally decreases as removal increases.

\begin{figure}[htbp]
    \centering
    \includegraphics[width=1\textwidth]{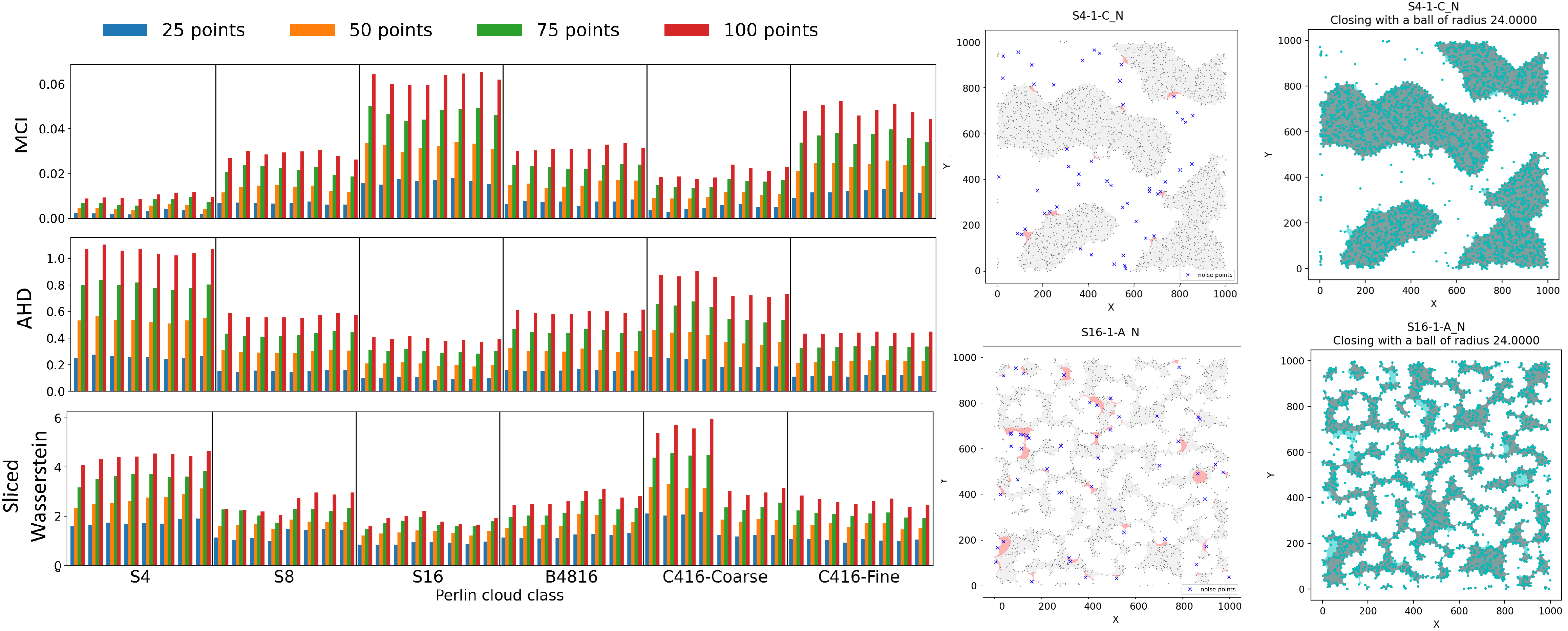}
    \caption{Response to impulsive noise: 25, 50, 75, or 100 points were added to the complementary white regions of each of the 48 Perlin clouds. Left: dissimilarity measures grouped by cloud class (HD95 is not displayed since got vale 0 in all the cases). Right: coarse (S4) and fine (S16) examples at closing radius \(r=24\); blue crosses mark added points and pink regions show newly reconstructed area.}
    \label{fig:noise}
\end{figure}

In the case of impulsive noise (see Fig. \ref{fig:noise}), added points that alter the reconstructed regions increase the MCI (though in very low quantities), whereas isolated points that leave the main structure largely unchanged have little impact and can more reasonably be regarded as noise. AHD shows the opposite tendency: it penalizes distant, isolated points more strongly, even when their structural effect is limited. MCI’s response is therefore more appropriate when the aim is to assess preservation of the underlying structure.

\paragraph{Implementation details and code and data availability.}
An implementation of MCI, as well as the full set of experiments and the needed files to test them, are provided in a Github repository\footnote{\url{https://github.com/afalvarez95/point-sets-structural-similarity}.}.
The dissimilarity measure MCI has been implemented using Python with the library Shapely \cite{Gillies2022Shapely}, which allows for morphological dilation and erosion operations. Since Shapely works with polygonal approximations of the closings, the areas of their unions and intersections, and hence the values of MCI, are approximations of the exact ones. The code allows for batch experimentation of a pair of point sets: either generated randomly, passed as an argument, or given by a CSV file.

The solution also includes an interface made with Gradio \cite{Abid2019Gradio}, which allows easy data visualization and comparison of the measurements used in the experiments. 

For the synthesis of the ground database we used the Python \texttt{perlin-noise} implementation of Perlin noise. We computed sliced Wasserstein distance using the Python Optimal
Transport (POT) library \cite{flamary2026pot,flamary2021pot}, with $p=2$ and
$L=50$ uniformly sampled projection directions with 20 random fixed seeds and then averaged.

\section{Conclusions}\label{sec:conclusions}

We introduced the Multiscale Closing Index (MCI), a dissimilarity between 2D point sets obtained by averaging, over a range of radii, the Jaccard distance between the morphological closings of the two sets. Its bounded range \([0,1]\) facilitates interpretation. It is symmetric, satisfies the triangle inequality, and converges to the Jaccard distance between the convex hulls as the maximal radius grows. In the experiments, compared with HD95, AHD and sliced Wasserstein, it separated samplings of the same underlying structure from samplings of different structures and more consistently reflected the impact of alterations relative to structural scale, generally increasing with severity while remaining tolerant of isolated noise. Adjustable minimum and maximum radii allow adaptation to different sampling densities and application-specific scales. On the other hand, this flexibility adds a parameter choice to the analysis pipeline. Another limitation is that it is not definite (distinct point sets may have zero index when restricted to a range of radii) for all radii ranges.

\begin{credits}
\subsubsection*{\ackname} This research was partially funded by the French National Research Agency (ANR) under the France 2030 program, grant
PR[AI]RIE-PSAI (ANR-23-IACL-0008). M.J.J. is supported by grant PID2024-155867NB-I00 funded by MICIU/AEI/10.13039/501100011033 and ERDF/EU and by the IMUS–María de Maeztu grant CEX2024-001517-M (Apoyo a Unidades de Excelencia María de Maeztu), funded by MICIU/AEI/10.13039/501100011033.

\subsubsection*{AI disclosure.}
In this work, we have used generative AI tools for improving the writing and readability, searching and completing the literature review and helping to summarize and fit to the 12 pages limit. The original idea of the method introduced in the paper has been developed completely on our own. Regarding implementation, though We have not used AI tools for generating the data, nor coding the main pipeline and experiments, Claude (Opus 5.5 medium effort) has been used for experiments replication testing, code debugging, readme generation and proofreading and requierements.txt compatible versioning. We have reviewed all AI-assisted work. We take responsibility for the final content of this work, including text, claims or artifacts produced with the aid of generative AI.
\end{credits}

%
%

 \bibliographystyle{splncs04}
 \bibliography{mybibliography}

\end{document}